\documentclass{article}

\usepackage[dvips]{graphicx}
\usepackage{amssymb}
\usepackage{amsmath}

\usepackage{caption}
\usepackage{hyperref}
\usepackage{arcs}
\usepackage{xcolor} 
\usepackage{subfigure}
\usepackage{authblk}
\usepackage{calrsfs}

\newtheorem{theorem}{Theorem}

\usepackage{graphicx} % Required for inserting images

\newenvironment{AMS}{\small\bf 2020 AMS subject classification: }{} 

\title{Unisolvence of Kansa collocation  for elliptic equations by  polyharmonic splines with random fictitious centers}

\author[1]{M. Mohammadi}
\affil[1]{Kharazmi University, Tehran}
\author[2]{A. Sommariva} 

\author[3]{M. Vianello\footnote{Corresponding author: marcov@math.unipd.it}}
\affil[1,2,3]{University of Padova, Italy}

\date{}

\begin{document}

\maketitle

\begin{abstract}
We make a further step in the unisolvence open problem for  unsymmetric Kansa collocation, proving nonsingularity of Kansa matrices with polyharmonic splines and random fictitious centers, for second-order elliptic equations with mixed boundary conditions. We also show some  numerical tests, where 
the fictitious centers are local random perturbations of predetermined collocation points.
\end{abstract}

\vskip0.5cm
\noindent
\begin{AMS}
{\rm 65D12, 65N35.}
\end{AMS}
\vskip0.2cm
\noindent
{\small{\bf Keywords:} second-order elliptic equations, mixed boundary conditions, meshless methods, unsymmetric Kansa collocation method, Radial Basis Functions, random fictitious centers, polyharmonic splines, unisolvence.}
\vskip1.5cm

\section{Introduction}

Strong form meshless collocation by Radial Basis Functions, named after the pioneering work of E.J. Kansa in the late 80's \cite{K86,K90}, has  been widely used in  numerical modelling of scientific and engineering problems for almost 40 years; cf., with no pretence of exhaustivity,  \cite{CDR20,CFC14,CKD19,CKA21,CL20,F07,SW06} and the references therein. Despite the  manifest success of Kansa method and several studies devoted to its theoretical and computational features, for example on unisolvence in the overtesting framework via least squares (cf. e.g. \cite{CLS18,HS01,LOS06,S16,SW06}), the basic problem of invertibility of unsymmetric square collocation matrices remains still substantially open. It is indeed well known since the fundamental paper by Hon and Schaback \cite{HS01} that there exist singular collocation designs (though they are extremely ``rare''), but finding sufficient conditions that ensure invertibility is a challenging topic. This theoretical lack was clearly recognized in the popular textbook \cite{F07}, and still only quite partially fixed. 

To be precise, some specific  results have been recently obtained for random collocation, in the standard (but quite  difficult) situation where RBF centers and collocation points coincide, namely in the particular case of the Poisson equation with Dirichlet boundary  conditions, cf. \cite{CDRDASV24,DASV24}. Nevertheless, while the assumptions in \cite{CDRDASV24}  (which treats the popular MultiQuadrics) are weak, in particular the boundary has no specific property besides those required for well-posedness of the differential problem, in \cite{DASV24}, dealing with Thin-Plate Splines (TPS) with no polynomial addition, the boundary is assumed to possess an analytic parametrization. 
All these restrictions make the quoted results not completely satisfactory, even in the framework of random collocation. 

In this paper, expanding  the proving technique presented in \cite{DASV23}, we focus on almost sure unisolvence of {\em random meshless collocation} by polyharmonic splines, trying to work with more generality on {\em second-order elliptic equations} and {\em mixed Dirichlet-Neumann boundary conditions}. The main idea that allows the generalization is to work with distinct center and collocation points, in particular the collocation points will be fixed and the centers randomly chosen. 
While the use of a sort of ``fictitious'' centers is not new in the literature on the Kansa method (essentially in the least squares framework), 
cf. e.g. \cite{CKD19,CKA21,ZLHW22}, the novelty here is that we deal with square collocation matrices. Moreover, the centers are independent but not necessarily identically distributed random variables, in particular we can take the centers as local random perturbations of predetermined collocation points. On the other hand, in the implementation we can exploit an appealing computational feature of polyharmonic splines, that of being {\em scale independent}, cf. e.g. \cite{F07}.

The main theoretical result is proved in Section 2, while in Section 3 we present some numerical tests to show the applicability of this kind of Kansa collocation, which is theoretically unisolvent in a probabilistic sense.

\section{Second-order elliptic equations with mixed boundary conditions}

In this paper, we consider second-order elliptic equations with variable coefficients and mixed boundary conditions 
\begin{equation}\label{cdr}
\left\{
\begin{array}{l}
\mathcal{L}u(P)=\sum_{i,j=1}^d{c_{ij}(P)\partial^2_{x_ix_j}} u(P)+\langle \nabla u(P),\vec{b}(P)\rangle+\rho(P)u(P)=f(P)\;,\;P\in \Omega\;,
\\ \\
\mathcal{B}u(P)=\chi_{\Gamma_1}(P)u(P)+\chi_{\Gamma_2}(P)\partial_\nu u(P)=g(P)\;,\;P\in \partial \Omega\;,
\end{array} 
\right .
\end{equation}
where $\Omega\subset \mathbb{R}^d$ is a bounded domain (connected open set), $P=(x_1,\dots,x_d)$, and the differential operator is {\em elliptic}, i.e. 
\begin{equation} \label{elliptic}
\sum_{i,j}{c_{ij}(P)\xi_i\xi_j}\neq 0\;,\;\;\forall P\in \Omega\;,\;\forall\xi\in \mathbb{R}^d\setminus \{0\}\;.
\end{equation}

\noindent Moreover,  
$\nabla=(\partial_{x_1},\dots,\partial_{x_d})$ denotes the gradient and $\langle \cdot,\cdot\rangle$ the inner product in $\mathbb{R}^d$, $\vec{b}$ is a vector field, 
$\partial_\nu=\langle \nabla,\vec{\nu}\rangle$ is the normal derivative at a boundary point, and $\chi$ denotes the characteristic function. The boundary is indeed splitted in two disjoint portions, namely $\partial \Omega=\Gamma_1\cup \Gamma_2$. If $\Gamma_2=\emptyset$
or $\Gamma_1=\emptyset$ we recover purely Dirichlet or purely Neumann conditions, respectively. 

We recall that polyharmonic splines correspond to the radial functions $$\phi(r)=r^{k} \log(r)\;,\;\;4\leq k\in 2\mathbb{N}$$ (TPS, Thin-Plate Splines, order $m=k/2+1$) and $$\phi(r)=r^k\;,\;\;2<k\in \mathbb{N}\;,\;\;k\;\mbox{odd}$$ (RP, Radial Powers, order $m=\lceil k/2\rceil$), cf. e.g. \cite{F07}. The restriction on the exponents in the present context guarantees the existence of second derivatives.

In the sequel, we shall use the following notation:
\begin{equation} \label{phiA}
\phi_A(P)=\phi(\|P-A\|)
\end{equation}
where $A=(a_1,\dots,a_d)$ is the RBF center and $\|\cdot\|$ the Euclidean norm; notice that, for both TPS and RP,  $\phi_A(P)$ is a {\em real analytic function} of $P$ for fixed $A$ (and of $A$ for fixed $P$) in $\mathbb{R}^d$ up to $P=A$, due to analyticity of the univariate functions $\log(\cdot)$ and $\sqrt{\cdot}$ in $\mathbb{R}^+$.

Observe that taking derivatives with respect to the $P$ variable
\begin{equation} \label{L}
\partial^2_{x_ix_j}\phi_A(P)=\delta_{ij}\frac{\phi'(r)}{r}+\frac{(x_i-a_i)(x_j-a_j)}{r^2}\left(\phi''(r)-\frac{\phi'(r)}{r}\right)\;,
\end{equation}
$$
\nabla \phi_A(P)=
(P-A)\phi'(r)/r\;,\;\;
\partial_\nu\phi_A(P)=\langle P-A,\vec{\nu}(P)\rangle\,\phi'(r)/r\;,
$$
where $r=\|P-A\|$.
Notice that in the TPS case
$$
\phi'(r)/r={r^{k-2}(k\log(r)+1)}\;,\;\;
\phi''(r)-\phi'(r)/r=r^{k-2}(k(k-2)\log(r)+2k-2)\;,
$$
whereas in the RP case
$$
\phi'(r)/r=kr^{k-2}\;,\;\;
\phi''(r)-\phi'(r)/r=k(k-2)r^{k-2}\;.
$$
%}
The Kansa collocation matrix can be written as 

$$
K_N=
\left(\begin{array} {ccccc}
\mathcal{L} \phi_{A_1}(P_1) & \cdots &  \cdots & \cdots & \mathcal{L}\phi_{A_N}(P_1)
\\
\\
\vdots & \vdots & \vdots & \vdots & \vdots
\\
\\
\mathcal{L} \phi_{A_1}(P_{N_I}) & \cdots &  \cdots & \cdots & \mathcal{L}\phi_{A_N}(P_{N_I})
\\
\\
\mathcal{B} \phi_{A_1}(Q_1) & \cdots &  \cdots & \cdots & \mathcal{B}\phi_{A_N}(Q_1)
\\
\\
\vdots & \vdots & \vdots & \vdots & \vdots
\\
\\
\mathcal{B}\phi_{A_1}(Q_{N_B}) & \cdots &  \cdots & \cdots & \mathcal{B}\phi_{A_N}(Q_{N_B})

\end{array} \right)
$$
where $N=N_I+N_B$, with $\{P_1,\dots,P_{N_I}\}$ distinct internal collocation points and $\{Q_1,\dots,Q_{N_B}\}$ distinct boundary collocation points.

We state now the main result. 
\begin{theorem}
Let $K_N$ be the polyharmonic spline Kansa collocation matrix defined above for the second-order elliptic equation (\ref{cdr}) with mixed boundary conditions on $\Omega\subset \mathbb{R}^ d$, $d\geq 2$, where $\{P_h\}\subset \Omega$ and $\{Q_k\}\subset \partial\Omega$ are any two fixed sets of distinct collocation points, and $\{A_i\}$ is a sequence of independent random points with possibly distinct probability densities $\sigma_i \in L^1_+(\mathbb{R}^d)$.

Then for every $N=N_I+N_B$ 
with $N_I\geq 0$, $N_B\geq 1$, the matrix $K_N$ is almost surely nonsingular.
\end{theorem}
\vskip0.2cm 
\noindent 
{\bf Proof.} 
For the induction base, let us take $N=1$, that is  $N_I=0$, $N_B=1$. 
Then, either 
$\det(K_1)=\phi_{A_1}(Q_1)$ when $Q_1\in \Gamma_1$, or $\det(K_1)=\partial_\nu \phi_{A_1}(Q_1)$ when $Q_1\in \Gamma_2$. 
In the first case $\phi_{A_1}(Q_1)=0$ with RP iff $A_1=Q_1$, an event which clearly has probability zero. With TPS we have also to consider that $A_1$ falls on the hypersphere $\|Q_1-A\|=1$ due to the logarithmic factor, again an event with probability zero.
In the second case, 
$\partial_\nu \phi_{A_1}(Q_1)=0$ in the RP instance iff $A_1=Q_1$, an event which clearly has probability zero, or $Q_1-A_1$ is orthogonal to $\vec{\nu}(Q_1)$, that is $A_1$ falls in the tangent space to the boundary point $Q_1$. Also this event has probability zero, since a tangent space is a null set in 
$\mathbb{R}^d$. In the TPS instance we have also to take into account the case $k\log(\|Q_1-A_1\|)=-1$, i.e. that $A_1$ falls on the hypersphere $\|Q_1-A\|=\exp(-1/k)$, again an event with probability zero. 

For the inductive step, in case an internal collocation point is added we consider the auxiliary matrix (as a function of the new random center $A$)
$$
K_I(A)=
\left(\begin{array} {cccccc}
\mathcal{L} \phi_{A_1}(P_1) & \cdots &  \cdots & \cdots & \mathcal{L}\phi_{A_N}(P_1) & \mathcal{L}\phi_{A}(P_1)
\\
\\
\vdots & \vdots & \vdots & \vdots & \vdots & \vdots
\\
\\
\mathcal{L} \phi_{A_1}(P_{N_I}) & \cdots &  \cdots & \cdots & \mathcal{L}\phi_{A_N}(P_{N_I}) & \mathcal{L}\phi_{A}(P_{N_I})
\\
\\
\mathcal{L} \phi_{A_1}(P_{N_I+1}) & \cdots &  \cdots & \cdots & \mathcal{L}\phi_{A_N}(P_{N_I+1}) & \mathcal{L}\phi_{A}(P_{N_I+1})
\\
\\
\mathcal{B} \phi_{A_1}(Q_1) & \cdots &  \cdots & \cdots & \mathcal{B}\phi_{A_N}(Q_1) & \mathcal{B}\phi_{A}(Q_1)
\\
\\
\vdots & \vdots & \vdots & \vdots & \vdots & \vdots
\\
\\
\mathcal{B}\phi_{A_1}(Q_{N_B}) & \cdots &  \cdots & \cdots & \mathcal{B}\phi_{A_N}(Q_{N_B}) & \mathcal{B}\phi_{A}(Q_{N_B})
\end{array} \right)
$$
with $K_{N+1}=K_I(A_{N+1})$, 
while in case a boundary  collocation point is added we consider the matrix
$$
K_B(A)=
\left(\begin{array} {cccccc}
\mathcal{L} \phi_{A_1}(P_1) & \cdots &  \cdots & \cdots & \mathcal{L}\phi_{A_N}(P_1) & \mathcal{L}\phi_{A}(P_1)
\\
\\
\vdots & \vdots & \vdots & \vdots & \vdots & \vdots
\\
\\
\mathcal{L} \phi_{A_1}(P_{N_I}) & \cdots &  \cdots & \cdots & \mathcal{L}\phi_{A_N}(P_{N_I}) & \mathcal{L}\phi_{A}(P_{N_I})
\\
\\
\mathcal{B} \phi_{A_1}(Q_1) & \cdots &  \cdots & \cdots & \mathcal{B}\phi_{A_N}(Q_1) & \mathcal{B}\phi_{A}(Q_1)
\\
\\
\vdots & \vdots & \vdots & \vdots & \vdots & \vdots
\\
\\
\mathcal{B} \phi_{A_1}(Q_{N_B}) & \cdots &  \cdots & \cdots & \mathcal{B}\phi_{A_N}(Q_{N_B}) & \mathcal{B}\phi_{A}(Q_{N_B})
\\
\\
\mathcal{B} \phi_{A_1}(Q_{N_B+1}) & \cdots &  \cdots & \cdots & \mathcal{B}\phi_{A_N}(Q_{N_B+1}) & \mathcal{B}\phi_{A}(Q_{N_B+1})
\end{array} \right)
$$
with $K_{N+1}=K_B(A_{N+1})$.
Now developing the determinants by the last column we get
$$
F(A)=\det(K_I(A))=\alpha_1\mathcal{L}\phi_A(P_1)+\dots+
\alpha_{N_I}\mathcal{L}\phi_A(P_{N_I})+\alpha_{N_I+1}
\mathcal{L}\phi_A(P_{N_I+1})
$$
$$
+\beta_1\mathcal{B}\phi_A(P_1)+\dots+
\beta_{N_B}\mathcal{B}\phi_A(P_{N_B})
$$
and 
$$
G(A)=\det(K_B(A))=\alpha_1\mathcal{L}\phi_A(P_1)+\dots+
\alpha_{N_I}\mathcal{L}\phi_A(P_{N_I})
$$
$$
+\beta_1\mathcal{B}\phi_A(Q_1)+\dots+
\beta_{N_B}\mathcal{B}\phi_A(Q_{N_B})
+\beta_{N_B+1}\mathcal{B}\phi_A(Q_{N_B+1})
$$
where $\alpha_h$, $\beta_k$ are the corresponding minors with the appropriate sign. Since $|\alpha_{N_I+1}|=|\beta_{N_B+1}|=|\det(K_N)|$ is almost surely nonzero by inductive hypothesis, in both cases the determinants $F(A)$ and $G(A)$ as functions of $A$ are almost surely not identically zero in $\mathbb{R}^d$, because the functions 
$$\{\mathcal{L}\phi_A(P_h),\mathcal{B}\phi_A(Q_k),\,1\leq h\leq N_I,\,1\leq k\leq N_B\}$$ are linearly independent in $\mathbb{R}^d$ (as functions of $A$). To prove linear independence, first we show that $\mathcal{L}\phi_A(P_h)$ is analytic in $\mathbb{R}^d$ up to the point $A=P_h$ where it has a singularity of some partial derivative, for both RP and TPS. Indeed, by (\ref{phiA})-(\ref{L}), 
we have that $\mathcal{L}\phi_A(P_h)=
\mathcal{L}_A^{\circleddash}\phi_A(P_h)$ where 
$$
\mathcal{L}_A^{\circleddash}\phi_A(P_h)
=\sum_{i,j=1}^d{c_{ij}(P_h)\partial^2_{a_ia_j}\phi_A(P_h)}-\sum_{j=1}^d{b_j(P_h)\partial _{a_j}\phi_A(P_h)}+\rho(P_h)\phi_A(P_h)
$$
is a {\em constant coefficient elliptic operator} (acting on the $A$ variables; notice that the direction of the vector $\vec{b}$ has been inverted). Now, if 
$f_h(A)=\mathcal{L}\phi_A(P_h)$ were analytic at $A=P_h$ it would be analytic in a neighborhood of $P_h$, where $\phi_A(P_h)$ would solve the elliptic equation (in the $A$ variables) $\mathcal{L}_A^{\circleddash}\phi_A(P_h)
=f_h(A)$. By a famous result in the theory of elliptic equations with analytic coefficients and data (cf. e.g. \cite[Thm.7.5.1]{H63}), then $\phi_A(P_h)$ would be analytic in such a neighborhood and thus at $A=P_h$, which is a contradiction. 

On the other hand,  
$\mathcal{B}\phi_A(Q_k)$ is either $\phi_A(Q_k)$ for $Q_k\in \Gamma_1$ or alternatively $\partial_\nu\phi_A(Q_k)$ for $Q_k\in \Gamma_2$, and thus is analytic everywhere up to $A=Q_k$. To prove this fact for $\partial_\nu\phi_A(Q_k)$, take simply its restriction to the line through the point $Q_k=(q_1,\dots,q_d)$ parallel to the $a_1$ axis, namely $a_s=q_s$, $s=2,\dots,d$. Then $\partial_\nu\phi_A(Q_k)=\nu_1(Q_k)(q_1-a_1)\phi'(r)/r$, from which we obtain   
$\partial_\nu\phi_A(Q_k)=\nu_1(Q_k)(q_1-a_1)|q_1-a_1|^{k-2}(k\log(|q_1-a_1|)+1)$ for TPS and $\partial_\nu\phi_A(Q_k)=\nu_1(Q_k)(q_1-a_1)|q_1-a_1|^{k-2}$ for RP. Both have clearly a singularity of the $(k-1)$-th derivative in the $a_1$ variable at $a_1=q_1$. 

Consequently, if $\{\mathcal{L}\phi_A(P_h), \mathcal{B}\phi_A(Q_k)\}$ were dependent (as functions of $A$), one of them with singularity
say in $C$ (where $C$ is one of the collocation points) could be written as linear combination of the others.  
But such linear combination would be analytic at $A=C$, since the collocation points $\{P_h\}$ and $\{Q_k\}$ are all distinct,  thus all the other
functions are analytic at $C$ and a linear combination of analytic functions
at a point is analytic at such point. So this is a contradiction.

Now, both $F$ and $G$ are analytic in the open connected set $U=\mathbb{R}^d\setminus (\{P_h\}\cup \{Q_k\})$, and are almost surely not identically zero there, otherwise by continuity they would be zero in the whole $\mathbb{R}^d$. By a basic theorem in the theory of analytic functions (cf. \cite{M20} for an elementary proof), their zero sets in $U$
have then null Lebesgue measure and thus null measure with respect to any probability measure with integrable density in $\mathbb{R}^d$. Since the collocation set is finite and hence trivially it is a null set, the zero sets of $F$ and $G$ in $\mathbb{R}^d$, say $Z(F)$ and $Z(G)$, also have (almost surely) null measure. Concerning the probability of the corresponding events and taking for example $F$ (the conclusion is the same for $G$), we can write that
$$ 
\mbox{prob}\{\det(K_{N+1})=0\}=\mbox{prob}\{F(A_{N+1})=0\}
$$
$$
=\mbox{prob}\{F\equiv 0\}
+\mbox{prob}\{F\not\equiv 0\;\&\;A_{N+1}\in Z(F)\}
=0+0=0\;,
$$ 
and the inductive step is completed. \hfill $\square$

\section{Numerical Results}
In this section, we implement unsymmetric Kansa collocation method for solving the
linear second-order elliptic PDE (\ref{cdr}) on $\Omega=(0,1)^2$, with both TPS and RP RBF.
In all the test problems, we fix the set of collocation points $Y=\{Y_j\}_{j=1,\ldots,N}$ corresponding to a uniform discretization tensorial grid of the square, in lexicographic order. The random fictitious centers $X$ are obtained by local random perturbation of the collocation points via additive uniformly distributed random points in $(-\delta,\delta)^2$, as 
\begin{eqnarray*}
    X=Y+(2* \, \mbox{rand}(N,1)-1)*\delta\;,
\end{eqnarray*}
using a MATLAB notation. 
The accuracy is measured by the average of root mean
square error (RMSE) 
\begin{eqnarray*}
\mbox{RMSE}=\sqrt{\frac{1}{N}\displaystyle\sum_{j=1}^{N}(u_j-\tilde{u}_j)^2}
\end{eqnarray*}
obtained by $m$ random centers arrays $\{X_l\}$, $l=1,\ldots,m$, where $u_j$ and $\tilde{u}_j$ are the exact and approximate
solutions at the collocation node {$Y_j$}, respectively. In the sequel, we consider $m=100$ trials.
%We also take $51$ grid points along each axis for plotting of figures
%%%%%%%%%%%%%%
\subsection{Test problem 1}
As first test problem, we consider  equation (\ref{cdr}) with $\mathcal{L}=\Delta$ (Poisson equation), over
the unit square $\Omega=(0,1)^2$, with purely Dirichlet boundary conditions that is $\Gamma_1=\partial \Omega$ and $\Gamma_2=\emptyset$, 
where $f$ and $g$ are defined by forcing the analytical solution to be $u(x_1,x_2)=\sin(2\pi x_1)+\cos(2\pi x_2)$. 

In Fig. \ref{utest1}, we plot on a finer grid the approximate concentration profile $\tilde{u}$, computed by Kansa method with $N=441$ and $\delta=0.01$.
In Tables \ref{t1} and \ref{t2} we report the average RMSE (Root Mean Square Errors)
of Kansa method with TPS and RP with a couple of exponents, and different values of the neighborhood radius $\delta$ (observe that $\delta=0$ corresponds to 
classical collocation with $X=Y$, that is of centers coinciding with the collocation points). An example of distribution of errors around the average RMSE is illustrated in Fig. 3. 

We can observe that for the smallest values of $\delta$ the errors have the same size of those corresponding to classical collocation with  $X=Y$ (which can be considered a limit case), whose unisolvence is however not covered by the present theory, that on the contrary ensures  almost sure unsolvence for any $\delta>0$. For larger values of $\delta$, e.g. $\delta=0.1$, the errors are not satisfactory, which could be ascribable to the fact that less centers fall near the boundary. 

\begin{figure}[!]
%  % Requires \usepackage{graphicx}
 \center{ \includegraphics[scale=0.35]{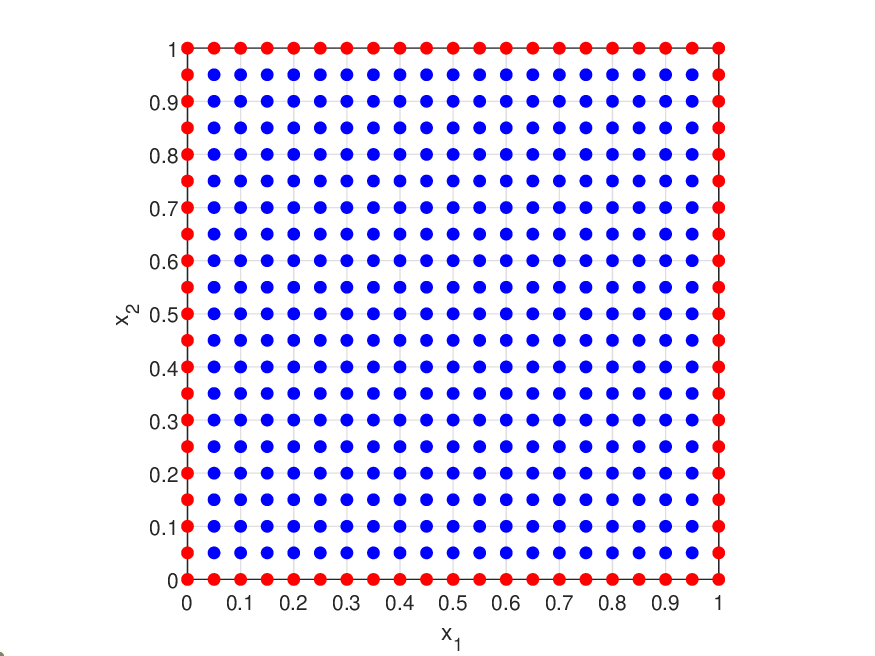} \includegraphics[scale=0.35]{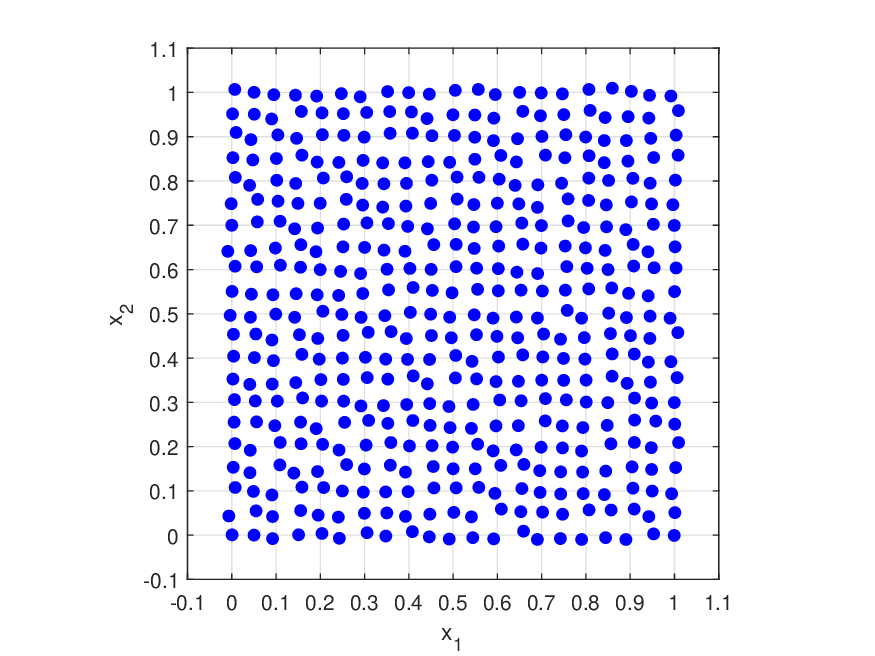}}
    \caption {$441$ collocation grid points (left) and the random fictitious centers distribution (right) for $\delta=0.01$ (Test problem 1).}\label{eqpointstest1}
    \end{figure}

    \begin{figure}[!]
%  % Requires \usepackage{graphicx}
 \center{\includegraphics[width=300pt]{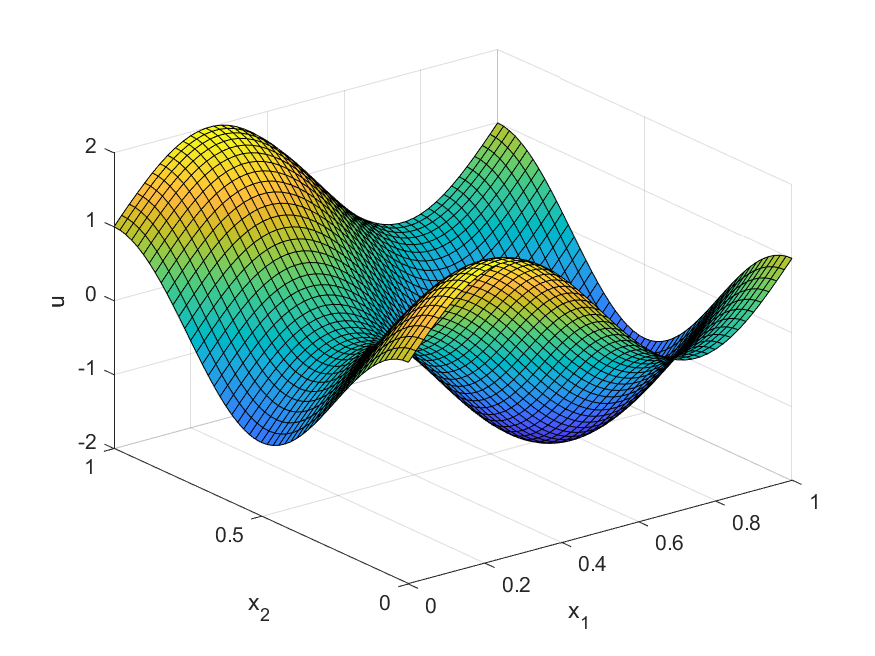}}
    \caption {Plot of the computed solution $u$ with $441$ collocation points and $\delta=0.01$ (Test problem 1).}\label{utest1}
    \end{figure}

\begin{table}[!]
  \centering  
  \caption{\footnotesize{Average of RMSE with Thin-Plate Splines and different values of $\delta$; test problem 1.}}  
  \label{t1}  
  \vspace{5pt} % Add vertical space before the table  
  \scriptsize{  
    \begin{tabular}{|c|c|c|c|c|||c|c|c|c|}  
      \hline  
      & \multicolumn{4}{c|||}{TPS $k=4$} & \multicolumn{4}{c|}{TPS $k=6$} \\
      \hline  
      $N$ & $\;\;\delta=0.1\;\;$ &  $\;\delta=0.01\;$ & $\delta=0.001$ & $\;\;\;\delta=0\;\;\;$ & $\;\;\delta=0.1\;\;$ &  $\;\delta=0.01\;$ & $\delta=0.001$ & $\;\;\;\delta=0\;\;\;$\\
      \hline  
      $121$ & 7.0e-01 & 9.8e-03 & 8.2e-03  & 8.2e-03& 4.4e-02 & 5.6e-03 & 4.6e-03 & 4.6e-03 \\
      \hline
      $441$  & 7.0e-01 & 1.4e-03 & 1.3e-03 & 1.3e-03 & 1.6e-02 & 3.7e-04 & 3.4e-04 & 3.4e-04 \\
      \hline
      $961$  & 7.6e+00 & 4.2e-04 & 4.0e-04 & 4.0e-04 & 8.7e-03 & 8.3e-05 & 7.3e-05 & 7.3e-05 \\
      \hline
      $1681$  & 1.2e+01 & 1.9e-04 & 1.8e-04 & 1.8e-04 & 4.6e-03 & 2.6e-05 & 2.4e-05 & 2.4e-05 \\
      \hline  
    \end{tabular}  
  }  
  \vspace{5pt} % Add vertical space after the table  
\end{table}

\begin{table}[!]
  \centering  
  \caption{\footnotesize{Average of RMSE with Radial Powers and different values of $\delta$; test problem 1.}}  
  \label{t2}  
  \vspace{5pt} % Add vertical space before the table  
  \scriptsize{  
    \begin{tabular}{|c|c|c|c|c|||c|c|c|c|}  
      \hline  
      & \multicolumn{4}{c|||}{RP $k=3$} & \multicolumn{4}{c|}{RP $k=5$} \\
      \hline  
      $N$ & $\;\;\delta=0.1\;\;$ &  $\;\delta=0.01\;$ & $\delta=0.001$ & $\;\;\;\delta=0\;\;\;$ & $\;\;\delta=0.1\;\;$ &  $\;\delta=0.01\;$ & $\delta=0.001$ & $\;\;\;\delta=0\;\;\;$\\
      \hline  
      $121$ & 1.0e+00 & 1.7e-02 & 1.5e-02 & 1.5e-02& 4.1e-01 & 5.8e-03 & 4.8e-03 & 4.8e-03 \\
      \hline
      $441$  & 1.2e+01 & 3.1e-03 & 3.2e-03  & 3.2e-03 & 7.5e-02 & 6.1e-04 & 5.5e-04 & 5.5e-04 \\
      \hline
      $961$  & 1.2e+01 & 1.2e-03 & 1.2e-03 & 1.2e-04 & 1.7e-01  & 1.7e-04 & 1.5e-04 & 1.5e-04 \\
      \hline
      $1681$  & 1.8e+01 & 5.5e-04 & 6.1e-04 & 6.2e-04 & 5.1e-02 & 6.8e-05 & 5.7e-05 & 5.7e-05 \\
      \hline  
    \end{tabular}  
  }  
  \vspace{5pt} % Add vertical space after the table  
\end{table}

\begin{figure}[!]
%  % Requires \usepackage{graphicx}
  \center{\includegraphics[scale=0.40]{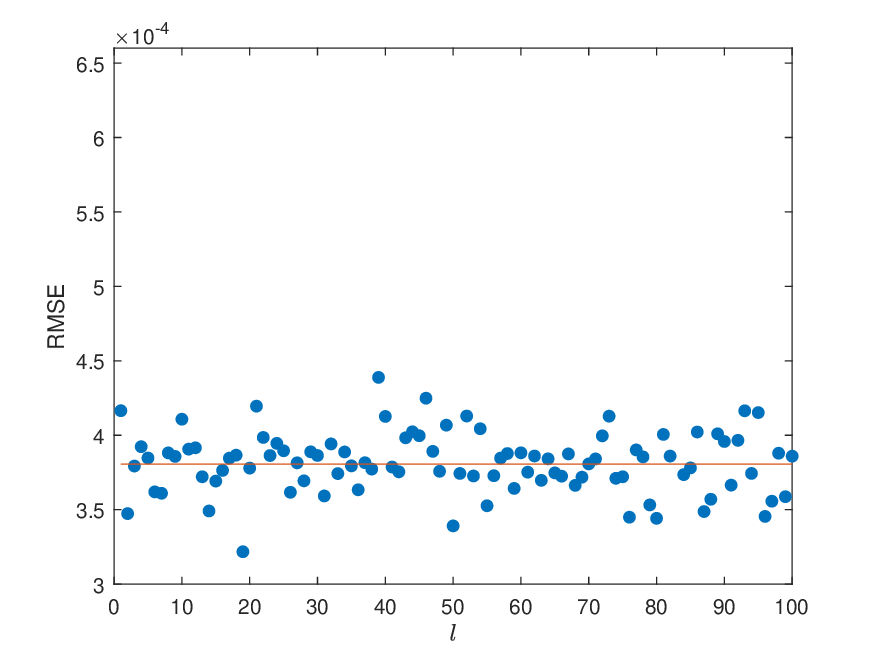}\includegraphics[scale=0.40]{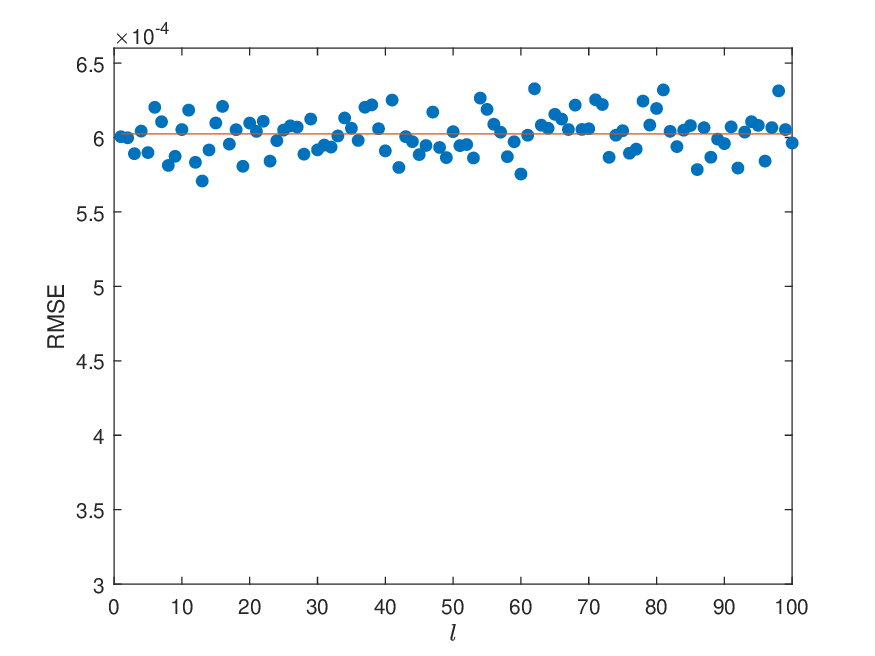}}
    \end{figure}
%%%%%%%%%%%%%%%%%
\begin{figure}[!]
%  % Requires \usepackage{graphicx}
  \center{\includegraphics[scale=0.40]{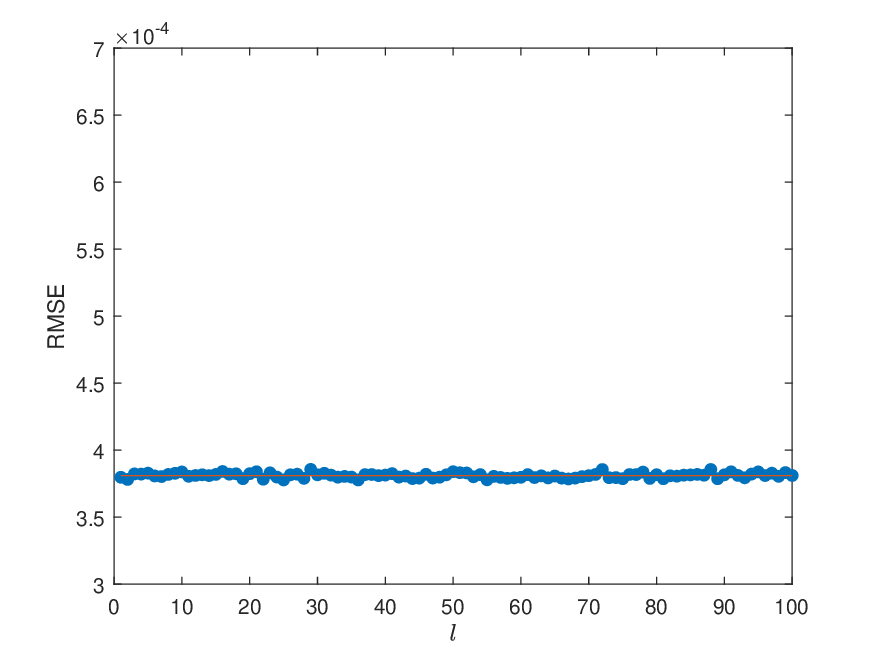}\includegraphics[scale=0.40]{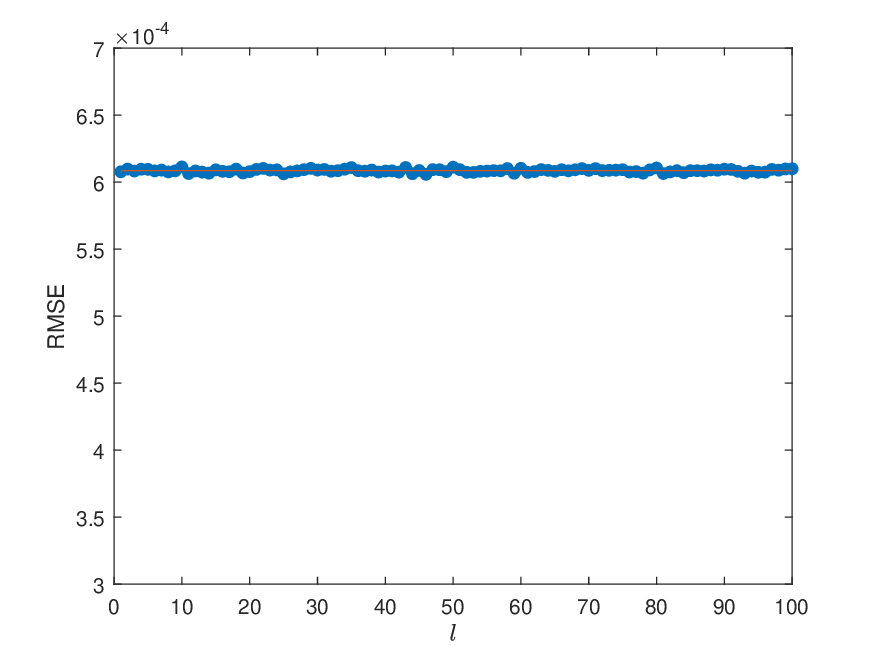}}
     \caption {Plot of  RMSE with $441$ equispaced  collocation points and  random centers arrays $X_l$, $l=1,\dots,100$,  for TPS RBF ($k=6$) (left) and RP RBF ($k=5$) (right), with $\delta=0.01$ (top) and $\delta=0.001$ (bottom); average RMSE in red, test problem 1.}\label{RMSETPSRP001}
    \end{figure}

\subsection{Test problem 2}
In the next experiment, we consider the convection-diffusion operator $\mathcal{L}u=\Delta u+\langle \nabla u,\vec{b}\rangle$, with  $\vec{b}=(1,1)$, over the unit square $\Omega=(0,1)^2$, and mixed boundary conditions in (\ref{cdr}), given by
\begin{eqnarray*}
    \Gamma_1&:&\{x_1=0,~ 0\leq x_2\leq 1\}\cup  \{x_1=1,~ 0\leq x_2\leq 1\},\\
    \Gamma_2&:&\{x_2=0,~ 0< x_1< 1\}\cup  \{x_2=1,~ 0< x_1< 1\}.
\end{eqnarray*}
Again, $f$ and $g$ are defined by selecting 
$u(x_1,x_2)=\sin(2\pi x_1)+\cos(2\pi x_2)$.

The numerical results corresponding to the application of the Kansa method with random fictitious centers are reported in Tables \ref{t3} and \ref{t4}, where we can see an error behavior similar to that discussed in Test problem 1. In particular, the smallest values of $\delta$ show errors having the same size of those corresponding to classical collocation with  $X=Y$. Almost sure unisolvence for any $\delta>0$ is again ensured by Theorem 1.

\begin{table}[!]
  \centering  
  \caption{\footnotesize{Average of RMSE with Thin-Plate Splines  and different values of $\delta$; test problem 2.}}  
  \label{t3}  
  \vspace{5pt} % Add vertical space before the table  
  \scriptsize{  
    \begin{tabular}{|c|c|c|c|c|||c|c|c|c|}  
      \hline  
      & \multicolumn{4}{c|||}{TPS $k=4$} & \multicolumn{4}{c|}{TPS $k=6$} \\
      \hline  
      $N$ & $\;\;\delta=0.1\;\;$ &  $\;\delta=0.01\;$ & $\delta=0.001$ & $\;\;\;\delta=0\;\;\;$ & $\;\;\delta=0.1\;\;$ &  $\;\delta=0.01\;$ & $\delta=0.001$ & $\;\;\;\delta=0\;\;\;$\\
      \hline  
      $121$ & 9.9e-01 & 1.2e-01 & 7.9e-02 & 7.2e-02& 1.1e+00 & 9.2e-02 & 8.7e-02 & 7.9e-02 \\
      \hline
      $441$  & 1.9e+00 & 7.1e-02 & 1.3e-02 & 1.3e-02 & 2.5e-01 & 1.3e-02 & 1.2e-02 & 1.1e-02 \\
      \hline
      $961$  & 3.3e+00 & 6.1e-02 & 5.1e-03 & 4.9e-03 & 1.5e-01 & 6.4e-03 & 3.6e-03 & 3.5e-03 \\
      \hline
      $1681$  & 5.0e+01 & 2.4e-02 & 2.7e-03 & 2.6e-03 & 7.6e-02 & 5.4e-03 & 1.5e-03 & 1.5e-03 \\
      \hline  
    \end{tabular}  
  }  
  \vspace{5pt} % Add vertical space after the table  
\end{table}

\begin{table}[!]
  \centering  
  \caption{\footnotesize{Average of RMSE with Radial Powers and different values of $\delta$; test problem 2.}}  
  \label{t4}  
  \vspace{5pt} % Add vertical space before the table  
  \scriptsize{  
    \begin{tabular}{|c|c|c|c|c|||c|c|c|c|}  
      \hline  
      & \multicolumn{4}{c|||}{RP $k=3$} & \multicolumn{4}{c|}{RP $k=5$} \\
      \hline  
      $N$ & $\;\;\delta=0.1\;\;$ &  $\;\delta=0.01\;$ & $\delta=0.001$ & $\;\;\;\delta=0\;\;\;$ & $\;\;\delta=0.1\;\;$ &  $\;\delta=0.01\;$ & $\delta=0.001$ & $\;\;\;\delta=0\;\;\;$\\
      \hline  
      $121$ & 3.8e+00 & 8.6e-01 & 4.5e-01 & 4.1e-01& 5.8e-01 & 3.7e-02 & 3.5e-02 & 3.2e-02 \\
      \hline
      $441$  & 1.7e+01 & 1.6e-01 & 1.1e-01 & 1.1e-01 & 6.5e-01 & 3.5e-02 & 9.5e-03 & 9.1e-03 \\
      \hline
      $961$  & 1.8e+01 & 1.2e-01 & 5.7e-02 & 5.5e-02 & 2.1e-01 & 6.7e-03 & 3.9e-03 & 3.8e-03 \\
      \hline
      $1681$  & 3.9e+01 & 1.1e-01 & 3.6e-02 & 3.5e-02 & 7.0e-01 & 3.3e-03 & 2.0e-03 & 2.0e-03 \\
      \hline  
    \end{tabular}  
  }  
  \vspace{5pt} % Add vertical space after the table  
\end{table}

\section{Conclusions}
We have proved almost sure unisolvence of unsymmetric Kansa collocation by polyharmonic splines with  random fictitious
centers, for second-order elliptic equations with mixed boundary conditions. Real analitycity up to the centers plays a key role. This result contributes to fill a theoretical gap in the framework of Kansa meshless collocation, since determining sufficient unisolvence conditions of square collocation systems by Radial Basis Functions is still a substantially open problem, despite the widespread adoption of the method for engineering and scientific problems in the last decades.

Since in the present approach the centers can be non identically distributed random variables, they can be chosen as local random perturbations of a set of fixed collocation points. The fact that the collocation points are fixed and distinct from the centers allows to manage general elliptic operators and boundary conditions, differently from some recent results on random collocation. For small perturbations, the numerical tests show an error behavior substantially equivalent to the classical case of centers coinciding with collocation points, where however sufficient conditions for unisolvence are not known. 
We also recall an appealing computational feature of polyharmonic splines, that of being {\em scale independent}. 

Extension to other RBF with finite 
smoothness, as well as to the popular case of Kansa collocation by MultiQuadrics which are infinitely differentiable everywhere (trying to exploit the presence of complex instead of real singularities), can be object of future research. Another interesting but challenging problem  concerns an extension of the theory to the case of quasi-random fictitious centers, namely to low-discrepancy sequences such as for example Halton points.

\vskip1cm 
\noindent
{\bf Acknowledgements.} 

Work partially supported by the DOR funds of the University of Padova, by the INdAM-GNCS 2024 Project “Kernel and polynomial methods for approximation and integration: theory and application software'' (A. Sommariva and M. Vianello) and by the INdAM-GNCS 2024 Project “Meshless Techniques and integro-differential equations: analysis and their application'' (M. Mohammadi). 

This research has been accomplished within the RITA ``Research ITalian network on Approximation", the SIMAI Activity Group ANA\&A, and the UMI Group TAA ``Approximation Theory and Applications".

\end{document}